%% file: GolodCodim2.tex
\documentclass[11pt]{amsart}
\usepackage[T1,T5]{fontenc}
\usepackage{geometry}                
\usepackage{graphicx}
\usepackage{amssymb}
\usepackage{epstopdf}
\input{preamble}

\usepackage{tikz-cd}
\usepackage{hyperref}

\usepackage{float}

\usepackage[all]{xy}

\def\phi{{\varphi}}
\def\gm{{\mathfrak m}}
\def\gn{{\mathfrak n}}
\DeclareMathOperator\codepth{{\rm codepth}}
\DeclareMathOperator\edim{{\rm embdim}}
\newtheorem{criterion}[theorem]{Criterion}

\title{Curves in \texorpdfstring{$\PP^{3}$}{P3} are Complete Intersections or Golod}
\date{October 31, 2026}                                           

\author[D. Eisenbud]{David Eisenbud}
\address[David Eisenbud]{Department of Mathematics, University of California, Berkeley, CA 94720}
\email{de@berkeley.edu}
\urladdr{eisenbud.github.io}

\begin{document}
\begin{abstract}
Let $I\subset S$ be an unmixed, codimension 2 ideal in a regular local ring, or an unmixed, codimension 2 ideal that is homogeneous in a positively graded polynomial ring over a field. We prove the statement of the title when $\pd_{S}(S/I) \leq 3$, 
and conjecture that this dichotomy holds in codimension 2 
without the restriction to projective dimension 3.
\end{abstract}

\maketitle

\section*{Introduction}
In this paper we will prove: 

\begin{theorem}\label{main}
If $S$ is a regular local ring and
$I\subset S$ is an unmixed ideal of codimension $\leq 2$
that is not a codimension 2 complete intersection, and the projective dimension of $S/I$ is $\leq 3$, then $S/I$ is Golod. 
The same is true if $S$ is a positively graded polynomial ring over a field and $I$ is a homogeneous ideal. \end{theorem}
  
Perhaps the most interesting new family of examples comes from curves in $\PP^3$:

\begin{corollary}\label{GolodCurves}
 Let $X$ be a curve in $\PP^3$ with homogeneous coordinate ring $A(X)$.
 If $X$ is not a complete intersection, then $A(X)$ is a Golod ring. \qed
\end{corollary}

Here we may take the definition of ``curve'' to be very broad: we require only a purely 1-dimensional subscheme (that is, with no embedded points) and we take $A(X)$ to be defined by the saturated homogeneous ideal of $X$. Since the ideal is saturated, $A(X)$ has depth $\geq 1$ and thus, by the Auslander-Buchsbaum formula,
projective dimension $\leq 3$, so the Corollary follows from Theorem~\ref{main}.

Although the Golod property has been studied intensively for some time, most of the results have been combinatorial; practically the only geometrically interesting families of examples known until now have been ideals of maximal minors of various families of linear matrices (see~\cite{BCD} and the many references there) and  codimension 2 lattice ideals \cite[Theorem 2.3 and Lemma 6.2]{PeevaSturmfels}.

Corollary~\ref{GolodCurves} suggests a  richer world of geometric examples: every nonsingular curve over an algebraically closed field can be embedded in $\PP^3$.
Furthermore, by \cite[Corollary~3.8 and Lemma~2.14(1)]{BCD}, unmixed monomial ideals of codimension 2 of any projective dimension that are not complete intersections are Golod. In light of our Theorem~\ref{main}, we conjecture:

\begin{conjecture}[Codimension 2 Dichotomy]
 Let $S$ be a regular local ring, and let $I\subset S$ be an unmixed ideal of codimension 2. If $I$ is not a complete intersection, then $S/I$ is Golod, and similarly when $S$ is a positively graded polynomial ring over a field and $I$ is homogeneous.
\end{conjecture}

Since the Golod property is stable under factoring out regular parameters (or linear forms) that are nonzerodivisors, it would suffice to prove the conjecture for the cases where $\dim S = d$ and $\pd_{S}(S/I) = d-1$, for each $d\geq 5$ (for example, surfaces in $\PP^{4}$, 3-folds in $\PP^{5}$\dots.)

\noindent\textbf{Outline of the proof of Theorem~\ref{main}:} Set $R=S/I$ and, in the local case, write $(S,\gm,k)$. Assume that $I$ is not a codimension 2 complete intersection. It was already known that if $I$ has codimension 1 or $I$ is perfect of codimension 2 then $S/I$ is Golod; the results and references are reviewed in Section~\ref{sec:start}. The new situation is where $I$ has codimension 2 and 
$\pd_S R = 3$. 

We will use Criterion~\ref{criterion}. This
applies only when $I\subset \gm^{2}$; but if $I$ contains a regular
parameter, then factoring it out we get an unmixed ideal of codimension 1, and such an ideal is principal, so $I$ was a complete intersection. 

Thus it suffices to show that the positive degree part $\Tor^S_{>0}(R,k)$ of the strictly graded-commutative algebra $\Tor^S_*(R,k)$ (with its homological grading)
 squares to 0. Starting from a non-vanishing product, we construct an associated prime of $I$ of codimension 3. The product $\Tor_1^S(R,k) \otimes \Tor_1^S(R,k) \to \Tor_2^S(R,k)$ is treated in Section~\ref{sec:first-mult}, while the product $\Tor_1^S(R,k) \otimes \Tor_2^S(R,k) \to \Tor_3^S(R,k)$ is treated in Section~\ref{sec:second-mult}. The codimension
bound in each case finally reduces to a generalization of Krull's Principal Ideal Theorem from \cite{EisenbudEvans},
described in Section~\ref{sec:tools}. \qed

One can also formulate Theorem~\ref{main} for rings that are not necessarily homomorphic images of regular local rings:
For a minimal regular Cohen presentation $Q\twoheadrightarrow\widehat R$, the Auslander-Buchsbaum formula gives
\[
\codepth R := \edim R-\depth R = \pd_Q(\widehat R).
\]
Thus $R$ is Golod if $\edim R-\dim R\leq 2$, $\codepth\, R\leq 3$, and $\widehat R$ is unmixed but not a complete intersection of codimension 2.
 
The possible Tor algebras of local rings of codepth $\leq 3$ were classified by \cite{WeymanClassification} and subsequently (despite the publication dates) in a different way, in~\cite{AKMClassification}. Golod rings of codepth 3 are of type $H(0,0)$ in Avramov's notation.  Unfortunately, the proofs shed no light on how geometric conditions such as the codimension might limit the range of possible Tor algebras.  Avramov used the result to prove the rationality 
of the Poincar\'e and Bass series for such rings in~\cite{AKMClassification} and \cite{AvramovCohomological}, and they form the beginning of the construction of generic Cohen-Macaulay algebras of codimension 3 by Weyman and his collaborators (see for example~\cite{GNW}).

  \textbf{Convention:} In this paper all rings will be Noetherian, and all modules considered will be finitely generated.
  
 \section{What is a Golod ring, and why the property is interesting}

\subsection{The definition of the Golod property} Golod rings were discovered by \cite{Golod} in connection with the Kaplansky-Serre problem, which asked about a  finiteness property of infinite resolutions. The standard definition goes by way of the Poincar\'e series, and we include it for the reader's convenience, though we will not use it.
For any local ring $(R,\gn,k)$ and finitely generated $R$-module $M$ we define the
$p$-th \emph{Betti number} of $M$ to be
$b^{R}_{p}(M) := \dim_{k}(\Tor_{p}^{R}(k,M))$
and define the ``Poincar\'e''  series of $M$ to be
$$
P^{R}_{M}(t) := \sum_{p=0}^{\infty}b_{p}^{R}(M)t^{p}.
$$
 Serre proved that if $R = S/I$ with $S$ regular local of dimension $e = \edim R$,
 then
$$
P^{R}_{k}(t)\preceq \frac{(1+t)^{e}}{1+t- tP^{S}_{R}(t)}
$$
where $\preceq$ denotes coefficient-wise inequality when the right hand side is expanded as a formal power series.
\begin{definition}\label{def:Golod}The local ring $R = S/I$ is Golod if coefficient-wise equality holds in the formula above.
\end{definition}
Avramov \cite[Proposition 6.5]{AvramovHopf} has shown that it suffices to check the equality of Betti numbers for
any single $p$ above $\codepth R := \edim R - \depth R$, so this very infinite definition leads to a finite criterion.

\def\BB{{\mathbb B}}
A recent result of Jesse Burke has put the Golod property into a different context, more relevant to this paper: If $R = S/I$, with
$(S,\gm,k)$ a regular local ring and $I\subset \gm^2$, Burke constructed an $R$-free resolution $\BB$ of $k$ such that the generating function for the ranks of free modules in $\BB$ satisfies the equality in the formula above, turning the numerical equality into the statement that $\BB$ is minimal. Moreover, the differentials in $\BB$ are constructed from the finite set of
maps that make up the $A_\infty$-algebra structure on a minimal $S$-free resolution of $R$. 

\begin{theorem} (\cite[Definition 6.1 and Theorem~6.13]{Burke})\label{Burke criterion}
 If $(S,\gm, k)$ is a regular local ring then a factor ring $R= S/I$, with $I\subset \gm^{2}$, is Golod if and only if all the
 $A_{\infty}$ structure maps on the truncation $\FF_{\geq 1}$ of the minimal $S$-free resolution $\FF$ of $R$  have image in $\gm \FF$.\qed
\end{theorem}

When $\pd_S(R)\leq 3$, as in Theorem~\ref{main}, the $S$-free resolution $\FF$ of $R$ can be given the structure of a differential graded algebra \cite{Gor3}
and the only $A_\infty$ maps are the differential and the multiplication maps.
 This leads to the characterization we will use:

\begin{criterion}\label{criterion}
If $(S,\gm,k)$ is a regular local ring and the projective dimension of $R = S/I$ is $\leq 3$,
with $I\subset \gm^{2}$,
then
 $R = S/I$  is Golod if and only if the multiplication maps
 $$\begin{aligned}
 (1)\quad \Tor_{1}^{S}(R,k) \otimes \Tor_{1}^{S}(R,k) &\to \Tor_{2}^{S}(R,k) \\
 (2)\quad \Tor_{1}^{S}(R,k) \otimes \Tor_{2}^{S}(R,k) &\to \Tor_{3}^{S}(R,k)
\end{aligned}
$$
are both zero.\qed
\end{criterion}
The sufficiency of the conditions is given in~\cite[Proposition 5.2.4(4)]{AvramovInfinite}, while Theorem~\ref{Burke criterion} gives both implications.

\subsection{The interest of the Golod property}

In the situation of Theorem~\ref{main}, consider the resolution of the residue field $k$ over 
$R:=A(X)$:
$$
\begin{tikzcd}
 \cdots\ar[r]& G_m\ar[r,"\psi_m"] &G_{m-1} \ar[r, "\psi_{m-1}"]&\cdots \ar[r] &G_1\ar[r, "\psi_1"] &R\ar[r]&k \ar[r] &0
\end{tikzcd}
$$
The \emph{$m$-th syzygy} of $k$ is by definition the image of $\psi_m$.

Recent work has shown that the Golod property of a local ring $R$ implies a strong restriction on the syzygies of the residue field of $R$.
 If $R$ is Golod, then \cite[Theorem 1.3]{DaoEtAl} shows that there are only finitely many 
indecomposable summands of the syzygy modules of the residue field, and these already occur as summands of one of 
$$
k = \syz_{0}^{R}(k),\gm_{R} = \syz_{1}^{R}(k), \syz_{2}^{R}(k), \dots, \syz_{e}^{R}(k)
$$
where $e = \edim R := \dim_{k}\gm_{R}/\gm_{R}^{2}$. In fact,
for every $i>e,$
$$
\syz^R_i(k)=\bigoplus_{1 \leq s \leq \pd_{S}R} \syz^R_{i-s-1}(k) ^{\rank F_s}\, , 
$$
where $F_{s}$ is the $s$-th free module in a minimal $S$-free resolution of $R$.
For example, if $R = S/(g)$ is a hypersurface ring, then the minimal free resolution of any module is eventually periodic of period 2~\cite{Eisenbud1980}, giving a still stronger restriction.

It would be interesting to know more about these indecomposable summands. In the graded setting $R = A(X)$ is the affine ring of a projective variety $X$. If $X$ is smooth,  the high syzygies of $k$ sheafify to a distinguished finite (up to twists) class of vector bundles on $X$. If $X$ is arithmetically Cohen-Macaulay, then the high syzygies are 
maximal Cohen-Macaulay $A(X)$-modules. 

\begin{example}
Let $S= \ZZ/101[w,x,y,z]$. Here are tables showing the summands of each of the first 5 syzygies of the residue field over some factor rings of projective dimension 1,2,3. The $i$-th row
of each table (counting from 0) represents the $i$-th syzygy; the columns represent the distinct
isomorphism classes, up to twist, of indecomposable
summands, and the $(i,j)$ entry in the matrix is the number of times the $j$-th
indecomposable module appears in the decomposition of $\syz_i^R(k)$. Subsequent rows of each table have the same number of columns and are determined by the formula above. \smallskip

\noindent \textbf{Codimension 1, The nonsingular quadric surface:}
 Let $R = S/(xy-wz)$, so that $X\cong \PP^{1}\times \PP^{1}$ is the smooth quadric in $\PP^{3}$. It is well-known that there are just two indecomposable non-free maximal Cohen-Macaulay $R$-modules $M_{1}, M_{2}$ up to twists \cite{Knoerrer}.
 
 $\left(\!\begin{array}{ccccc}
      1&0&0&0&0\\
      0&1&0&0&0\\
      0&0&1&0&0\\
      0&0&0&2&2\\
      0&0&0&2&2
      \end{array}\!\right)$.
      
For example, the last row of the table indicates that the 4th $R$-syzygy of $k$ is the direct sum of two copies of each of $M_{1}$ and $M_{2}$.
 
\smallskip

\noindent\textbf{Codimension 2, Cohen-Macaulay -- the twisted cubic:}

$\PP^1 \ni [s,t] \mapsto [s^3, s^2t, st^2, t^3]\in \PP^3$

$\left(\!\begin{array}{ccc}
       1&0&0\\
       0&1&0\\
       0&0&3\\
       0&0&6\\
       0&0&12
       \end{array}\!\right)$
\smallskip

\noindent\textbf{Codimension 2,  not Cohen-Macaulay, The smooth rational quartic:}

$\PP^1 \ni [s,t] \mapsto [s^4, s^3t, st^3, t^4]\in \PP^3$

$\left(\!\begin{array}{ccccc}
       1&0&0&0&0\\
       0&1&0&0&0\\
       0&0&1&0&0\\
       0&3&0&2&2\\
       0&4&3&2&2
       \end{array}\!\right)$

\end{example}

Here is another application: 
If the curve is reduced, irreducible and nondegenerate then, 
over an algebraically closed field of characteristic 0, its general hyperplane section is a set of points in linearly general position, but the ideal-theoretic  intersection will typically have an embedded ``irrelevant'' component as well as the reduced points. Corollary~\ref{GolodCurves} gives an obstruction to a certain type of deformation of the homogeneous coordinate ring of such a set of points:

\begin{corollary}
Let $S = k[x,y,z,w]$, $S' = k[x,y,z]$,  and  let $I' \subset S'$ be a homogeneous codimension 2 ideal (possibly with $(x,y,z)$ as associated prime) that is not a complete intersection. Set $R' = S'/I'$. If $R'$ can be written as $R/(\ell)$ for
some ring $R$ of the form $S/I$ with $I$ unmixed of codimension 2
such that $\ell$ is a linear nonzerodivisor on $R$,  then $R'$ is Golod.
\end{corollary}
\begin{proof}
By Corollary~\ref{GolodCurves} the ring $R$ would be Golod, and this property is preserved on factoring out a linear nonzerodivisor \cite[Proposition~5.2.4(2),(3), p.~46]{AvramovInfinite}.
\end{proof}

It is interesting to compare our result to a codimension 3 result of
Christensen, Veliche, and Weyman
\cite[Theorem~1.2]{CVW}. Let $S$ be a regular local ring and let
$J\subset S$ satisfy
\[
\height J=\pd_S(S/J)=3,
\]
that is, $J$ is perfect of grade three. They show that a finite sequence
of complete-intersection links leads from $J$ to an ideal $J'$ such that
$S/J'$ is either a complete intersection or Golod. 

\subsection{Origin of the paper and AI contribution:}\label{sec:ack}

On the eve of a lecture on infinite free resolutions  that I was to give
 to an audience of algebraic geometers, I wondered what geometrically interesting classes of algebraic varieties might be Golod. I began by
 trying some easy classes of curves. To my surprise, every example I could construct that was not a complete
 intersection had the Golod property, and I was eventually able to prove the first part of Criterion~\ref{criterion}. However, I was stuck on the proof of the second part of the criterion until I explained
 the problem to ChatGPT 6.1, and it found the trick that is used at the end of this paper. I have simplified the AI contributions, and the exposition of the proofs is entirely mine.
 
\subsection{Acknowledgements: }\label{ack2} I am grateful to Kevin Tucker for helping me use AI to search for counterexamples; the failure
 of this search started me on the path to this paper. Xianglong Ni was
helpful in another attempt to find a counterexample. I'm grateful to Mohammed Abouzaid for encouraging me to
 trust AI -- which pushed me over the edge to getting AI help in proving the second part of the criterion.  
 
 \section{Previous results and reduction to the main case of Theorem~\ref{main}}\label{sec:start}
  
\noindent\textbf{Proof of Theorem~\ref{main}:}
 
 If $I\subset S$ has codimension 0, with $S$ regular, then $I=0$, and $S/I = S$ is Golod.
 
 If $I$ has codimension 1 then  $R$ is Golod by~\cite[Corollary 2 of Section 5]{Shamash}; see~\cite[Proposition 5.2.5]{AvramovInfinite} for a simple proof. Henceforward we assume
that $I$ has codimension 2, and thus $\pd_{S}(R) \geq 2$, but $I$ is not a complete intersection. 

Let $R = S/I$ be as in Theorem~\ref{main} with $\codim I = 2$. For simplicity we will treat the case when $S$ and $R$ are local; the proofs in the graded polynomial case are essentially identical.

 We will write the minimal $S$-free resolution of $R$ as
  \begin{center}
\begin{tikzcd}
\FF:\quad & 0\ar[r]  & F_{3}\ar[r,"\phi_{3}"] & F_{2} \ar[r,"\phi_{2}"] &F_{1} \ar[r,"\phi_{1}"] &S \ar[r]& R \ar[r] &0.
 \end{tikzcd}
\end{center}
We allow the case $F_3=0$. 

We will consider $\FF$ as a DG-algebra. Multiplication by an element $e\in F_{1}$
defines a null-homotopy $\sigma_{e}: \FF \to \FF[1]$ for multiplication by $\phi_{1}(e) \in S$. Any two null-homotopies for $\phi_{1}(e)$ are themselves homotopic, so to prove
that the conditions of Criterion~\ref{criterion} hold, it suffices to show, for all $i\geq 1$, that the image
of a null-homotopy $F_{i}\to F_{i+1}$  lies in $\gm \FF$.

In the case $\pd_{S}(R) = 2$, the Hilbert-Burch Theorem shows that $I$ is generated by the $n\times n$ minors of an $n\times(n+1)$
matrix and one can give a formula discovered by Juergen Herzog in his 1973 Habilitationsschrift for the multiplication map, showing that $R$ is Golod if $n\geq 2$ \cite[5.3.4]{AvramovInfinite}. 

%
%
%

Henceforward we may assume that the projective dimension of $R$ is exactly 3. 
Also, as remarked in the introduction, $I\subset \gm^{2}$; for if $I$ contained
a regular parameter $g\in S$ then $I/(g)\subset S/(g)$ is unmixed of codimension 1, and 
since $S/(g)$ is regular, this implies that $I/(g)$ is principal, whence $I$ is a complete intersection.
Thus it suffices to check the two conditions of Criterion~\ref{criterion}.

\section{Tools}\label{sec:tools}
 
Before proving the two conditions of Criterion~\ref{criterion} in this case, we pause to remind the reader of four results we will use. Throughout this section, $S$ denotes a regular local ring,
though the results quoted are valid more generally.

The first two results we need are versions of
 the Generalized Principal Ideal Theorem from \cite{EisenbudEvans}\footnote{Proven in the cases we will use in \cite{EisenbudEvans} and in a more general case, with a simpler proof,  in \cite{Bruns} and \cite{EHU}.}. 

For an $R$-module $M$ and $x\in M$, we define the \emph{order ideal} $M^{*}(x)$ to be the ideal of $R$
generated by the images of $x$ under all homomorphisms $M \to R$. 
To extend the familiar concept of rank to an arbitrary Noetherian ring, we define
the \emph{rank} of $M$ to be the maximum, over the minimal primes $P$ of $R$,
 of the minimal number of generators of $M_P$ as an $R_P$-module. 

\begin{theorem}(Generalized Principal Ideal Theorem)\label{gpit1}
Let $x\in M$ be an element of a module  over a Noetherian ring $R$.
If there is a prime ideal $Q$ of $R$ such that $x \in Q M_{Q}$, then
$$
\codim M^{*}(x) \leq \rank(M),
$$
\qed
\end{theorem}

For example, if $(R,\gm)$ is local, $M = R^{n}$, and $x = (x_{1}, \dots x_{n})\in \gm R^{n}$, then
$M^{*}(x) = (x_{1}, \dots, x_{n})\subset R$, recovering Krull's version. The next result (and a similar statement for minors of any size) follows
by a short argument from Theorem~\ref{gpit1}:

\begin{theorem}\label{gpit2}
 Suppose that $\phi$ is a $p\times q$ matrix, with $p\geq q$, over a
 catenary, equidimensional local ring $S$  and $\phi'$ is the $(p-1)\times q$
 submatrix consisting of the first $p-1$ rows of $\phi$. If the last row of
 $\phi$ has entries in the maximal ideal of $S$, then
 $$
 \codim (I_{q}(\phi)) \leq 1+\codim (I_{q}(\phi')).
 $$
 (If $p=q$ then we take $I_{q}(\phi') = 0$.)
 \end{theorem}

These theorems were proven in the case we will use in \cite[Theorems 1.1 and 2.1] {EisenbudEvans}; a simpler and more general proof was subsequently given in \cite[Theorem 1 and Corollary 1] {Bruns}.

The next two results have to do with the maximal non-vanishing minors of the
maps in a free resolution. First, one implication from the criterion of exactness and the analysis of associated primes  from \cite{WMACE}. We make the convention that the determinant of a matrix of size 0 is 1.

\begin{theorem} \cite{WMACE}\label{WMACE}
Suppose that
 \begin{center}
\begin{tikzcd}
\FF:\quad & 0\ar[r]  & F_{n}\ar[r,"\phi_{n}"] & F_{n-1} \ar[r] & \cdots \ar[r] &F_{1} \ar[r,"\phi_{1}"] &S \ar[r]& S/I \ar[r] &0.
 \end{tikzcd}
\end{center}
is a finite $S$-free resolution. If $I(\phi_{i}) := I_{\rank \phi_{i}}(\phi_{i})$ denotes the ideal generated by the minors of $\phi_i$ of size $\rank(\phi_{i})$, then
\begin{enumerate}
\item $\sqrt{I(\phi_{1})} \subseteq  \sqrt{I(\phi_{2})} \subseteq  \cdots \subseteq \sqrt{I(\phi_{n})}$
\item $\codim (I(\phi_{i}))\geq i$
\item (\cite[Theorem~1.1(1) and Remark~(2), pp.~215--217]{EHV}) The associated primes of $I$  of codimension $c$ are precisely the primes
of codimension $c$ that contain $I(\phi_{c})$. In particular $I$ is unmixed
if and only if $\codim (I(\phi_{c})) > c$ for all $c>\codim I$.
\end{enumerate}
\end{theorem}

The last result we need is a special case of the ``First Structure Theorem'':
\begin{theorem}\cite[Theorem~3.1]{structure} \label{structure}
Suppose that $\FF$ is an $S$-free resolution of $S/I$ having length 3, as pictured above.
Set $r := \rank \phi_2$, so $r = \rank F_1 -1$. If
the image $I$ of $F_{1} \to S$ has grade $\geq 2$ then
$$
I\cdot I_{\rank \phi_{3}} (\phi_{3}) = I_r(\phi_{2}).
$$
More precisely, the $r\times r$ minor of $\phi_{2}$ omitting the $i$-th row
and columns with numbers in the set $J$ is equal, up to a unit, to the product of the $i$-th generator of $I$ and
the minor of $\phi_{3}$ made from rows in the set $J$.
\end{theorem}

\begin{example}\label{ex:structure}
Consider the ring
 $R = S/I$, where
$I = (w^2, xy^2, x^3, wx^2)\subset S = k[x,y,z,w],$ so $R$ is a
ring of projective dimension 3 and codimension 2. Only the three variables $w,x,y$ are involved, so  
the fact that $R$ has projective dimension 3 implies that $(w,x,y)$ is an associated prime of $I$,
which is therefore mixed (a more elementary reason: $(w,x,y)wxy\subset I$).
Here is the minimal free resolution of $R$ as an $S$-module:
\[
\resizebox{\linewidth}{!}{$\displaystyle
S^{2}\,\xrightarrow{\phi_3 =\begin{array}{l}\vphantom{y^{2}0}\\\vphantom{0y^{2}}\\\vphantom{-x0}\\\vphantom{w-x}\\\vphantom{0w}\end{array}\left(\!\begin{array}{cc}
      \vphantom{\left\{5\right\}}-x&0\\
      \vphantom{\left\{4\right\}}0&y^{2}\\
      \vphantom{\left\{4\right\}}y^{2}&0\\
      \vphantom{\left\{5\right\}}w&-x\\
      \vphantom{\left\{5\right\}}0&w
      \end{array}\!\right)}\,
      S^{5}\,\xrightarrow{\phi_2 =\begin{array}{l} \vphantom{-x^{2}0-x\,y^{2}00}\\ \vphantom{w-x0-y^{2}0}\\ \vphantom{0w00-y^{2}}\\ \vphantom{00w^{2}w\,xx^{2}}\end{array}\left(\!\begin{array}{ccccc}
     -x\,y^{2}&0&-x^{2}&0&0\\
      w^{2}&0&0&w\,x&x^{2}\\
      0&w&0&0&-y^{2}\\
      0&-x&w&-y^{2}&0
      \end{array}\!\right)}\,S^{4}
      \, \xrightarrow{\phi_1 =\left(
\!\begin{array}{cccc}
      w^{2}&x\,y^{2}&x^{3}&w\,x^{2}
      \end{array}  
      \!\right)}\,S^{1}
      $}
\]
 Note that the first two entries of $\phi_1$ form a regular sequence,
and the first column of $\phi_2$ is their Koszul syzygy. As an example of the structure theorem of \cite{structure},
the determinant $wx^3y^2$ of the matrix obtained by omitting the first two columns and the last row of $\phi_2$
is, up to sign, the product of the last entry, $wx^2$, of $\phi_1$ and the determinant, $-xy^{2}$, of the matrix formed by the first two
rows of $\phi_3$.
\end{example}


\section{\texorpdfstring{$\Tor_{1}\otimes \Tor_{1} \to \Tor_{2}$}{Tor1 tensor Tor1 to Tor2}}\label{sec:first-mult}

\begin{theorem}\label{thm:first-mult} Let $(S,\gm,k)$ be a regular local ring and let $I\subset \gm^{2}$ be an ideal of projective dimension $\leq 2$. Set $R= S/I$.
 If the product
 $$
\mu_{1,1}: \Tor_{1}^{S}(R,k) \otimes_{k} \Tor_{1}^{S}(R,k)  \to \Tor_{2}^{S}(R,k) 
$$ 
is non-zero, then either $I$ is a complete intersection or $I$ has a primary component of codimension 3.
\end{theorem}

The same proof works if $S$ is a positively graded polynomial ring, $\gm$ is the maximal
homogeneous ideal, and $I$ is homogeneous.

\begin{proof}
If $R$ has projective dimension 1 then the hypothesis cannot be satisfied. If $I$ is a complete intersection there is nothing to prove. 
If $R$ has 
projective dimension 2 then $R$ is Golod by~\cite[5.3.4]{AvramovInfinite} so $\mu_{1,1}=0$. Thus we may assume that
$R$ has projective dimension 3, and thus $\codim I\leq 3$. If $\codim\ I =1$ then $R$ is Golod by 
\cite[Proposition 5.2.5]{AvramovInfinite}, so the hypothesis contradicts Criterion~\ref{criterion}.  If $\codim\ I = 3$ then the conclusion is trivial.

Thus we may assume that $\codim I = 2$ and the minimal $S$-free resolution of $R$ has the form 
\begin{center}
\begin{tikzcd}
\FF:\quad & 0\ar[r]  & F_{3}\ar[r,"\phi_{3}"] & F_{2} \ar[r,"\phi_{2}"] &F_{1} \ar[r,"\phi_{1}"] &S \ar[r]& R \ar[r] &0.
 \end{tikzcd}
\end{center}

 As observed in \cite{Gor3},
 the resolution $\FF$ can be given the structure of a DG algebra, and this structure induces the algebra structure of $\Tor^{S}(R,k) = k\otimes_{S}\FF$.
 
Supposing that some product in the Tor algebra is nonzero, we may assume after a change of basis that  $e_{1} *e_{2} =  f\notin \gm F_{2}$, where $e_{1}, e_{2}$ are the first two elements in a basis of $F_{1}$ and $f$ is the first basis vector of $F_{2}$.  Since $\phi$ acts as a graded derivation on $\FF$,
 $$
 \phi_{2}(f) = \phi_{1}(e_{1})e_{2}- \phi_{1}(e_{2})e_{1}.
 $$
 Since $f$ is a minimal generator of $F_{2}$, this is among the minimal relations on $I$, so after reordering the columns of the matrix $\phi_1$ and the rows and columns of $\phi_{2}$
 we may assume that the two entries of the first column of a matrix for $\phi_2$ are $-\phi_{1}(e_{2}), \phi_{1}(e_{1})$ (see Example~\ref{ex:structure}).

The elements $g_{1}:=\phi_{1}(e_{1})$ and $g_{2}=\phi_{1}(e_{2})$ must be a regular sequence since otherwise one could divide the given syzygy by a non-unit. Since $I$ is not a complete intersection, we have $(g_{1}, g_{2})\subsetneq I$.
  
 Let $J \subsetneq (g_{1}, g_{2})$ be the ideal generated by the maximal non-vanishing minors of $\phi_{2}$ involving the first column. Let $\phi_{3}'$ be the matrix obtained from $\phi_{3}$
 by omitting the first row.  Since $\rank F_3 < \rank F_2$, the maximal minors of $\phi_3'$ are among the maximal minors of $\phi_3$.
 
 By the ``first structure theorem'' of \cite[Theorem 3.1]{structure}, $J$ is the product of $I$ and the ideal $K$ that is generated by the
 minors of $\phi_{3}'$ of size equal to $\rank F_{3}$. These are maximal minors, since
 $\rank F_{3}<\rank F_{2}$.

Since $(g_{1}, g_{2})\subsetneq I$ and $KI\subset (g_{1},g_{2})$, the ideal $K$ consists of zero divisors mod $(g_{1},g_{2})$.
Thus $K\subset P$ for some associated prime $P$ of  $(g_{1},g_{2})$. Since $(g_{1},g_{2})$
 is a complete intersection, it is unmixed, whence $\codim K \leq \codim P = 2$.
 
Since $\phi_{3}$ contains only one more row than $\phi_{3}'$, Theorem~\ref{gpit2} shows that the codimension of the ideal of maximal minors of $\phi_{3}$ is at most 1 more than that of $\phi_{3}'$; thus 
$\codim I(\phi_{3})\leq 3$. By Theorem~\ref{WMACE}(3), every codimension 3 minimal prime of $I(\phi_3)$ is an associated prime of $I$, as required.
\end{proof}

\section{\texorpdfstring{$\Tor_1\otimes\Tor_2\to \Tor_3$}{Tor1 tensor Tor2 to Tor3}}
\label{sec:second-mult}

\begin{theorem}
Let $(S, \gm, k)$ be a regular local ring and let $I\subset \gm^{2}$
be an ideal of codimension 2 such that $\pd_S(S/I) = 3$. If the multiplication map
 $$
\mu_{1,2}:  \Tor_1^{S}(S/I,k) \otimes\Tor_2^{S}(S/I,k) \to \Tor_3^{S}(S/I,k)
 $$
 is nonzero, then $I$ has an associated prime of codimension 3.
 \end{theorem}
 
\begin{proof}
With notation as before, let  $\FF, \phi$  be the minimal free resolution of $R$, regarded as 
a strictly graded-commutative DG-algebra.
The multiplicative structure on $\Tor^{S}_*(S/I,k)$ is induced from this DG algebra structure.
By hypothesis there exist $e\in F_1$ and $x\in F_2$ such that $y = ex \notin \gm F_3$.
Set $g= \phi_1(e) \in I$ and write $\sigma$ for multiplication by $e$ on $\FF$, so that 
$\sigma$ is a null-homotopy for multiplication by $g$; that is, $\phi \sigma + \sigma \phi$
is multiplication by $g$ on $\FF$. Note that because $e\in F_1$ has odd degree, $\sigma^2 = 0$.

Write $F_1 = G_1\oplus Se$. Since $Se$ maps isomorphically onto the ideal $Sg\subset S$, the image of $\phi_2$
does not meet $Se$, so the
minimal free resolution of $I/(g)$
may be written
\begin{center}
\begin{tikzcd}
0\arrow[r]& {F}_{3}\arrow[r,"\phi"] &{F}_{2}\arrow[r,"\phi' "]&{G}_{1}\arrow[r] &I/Sg \arrow[r] &0
\end{tikzcd}
\end{center}
where $\phi'$ is the composition of $\phi$ with the projection to $G_1$. Let $\sigma': G_1 \to F_2$
denote the composition of the null-homotopy $\sigma$ for $g$ and the inclusion of $G_1$. It is easy 
to check that $\phi'\sigma'$ is multiplication by $g$ on $G_1$.

Let
\[
T := \begin{pmatrix}
 \sigma\\\phi'
\end{pmatrix}:  F_{2}\to 
 F_{3} \oplus  G_{1},
\]
and set $N := \coker(T)$.
The map $\sigma \phi_3$ is multiplication by $g$ on $F_3$, so $\sigma$ is a monomorphism
on the image of $\phi_3$. Since the kernel of $\phi'$ is the image of $\phi_3$,
the map $T$ is a monomorphism.

Note that the image of $\sigma: S \to F_1$ is $Se$. If $b\in G_1\subset F_1$, then
$\phi_2 \sigma b +\sigma \phi_1 b = gb$ so $\phi'\sigma' b = gb$.
Also, $\sigma\sigma' = \sigma^2\mid_{G_1} =0$, so if $(a,b)\in F_3\oplus G_1$ then
$(ga,gb) = T(\phi_3(a)+\sigma'(b))$. Thus $N$ is annihilated by $g$.

Now consider the  diagram:
\begin{tiny}
\begin{center}
\begin{tikzcd}[ampersand replacement=\&, row sep=3.5em, column sep=3em]
  \&0\ar[d]\&0\ar[d]\&0\ar[d]\\
  0 \arrow[r]
    \& F_3 \arrow[r, "\phi_3"] \arrow[d, "g"]
    \& F_2 \arrow[r, "\phi'"]
      \arrow[d, "{\begin{pmatrix}\sigma\\\phi'\end{pmatrix}}"]
    \& \operatorname{im}(\phi') \arrow[d, hook] \arrow[r]
    \& 0 \\
  0 \arrow[r]
    \& F_3 \arrow[r, "{\begin{pmatrix}1\\0\end{pmatrix}}"] \arrow[d]
    \& F_3\oplus G_1 \arrow[r, "{\begin{pmatrix}0 & 1\end{pmatrix}}"] \arrow[d]
    \& G_1 \arrow[d] \arrow[r]
    \& 0 \\
  0 \arrow[r]
    \& F_3/gF_3 \arrow[r]\ar[d]
    \& N \arrow[r]\ar[d]
    \& I/(g) \arrow[r]\ar[d]
    \& 0 \\
    \&0 \&0 \&0
 \end{tikzcd}
\end{center}
\end{tiny}
    
The top two rows and all three vertical columns are exact, so the third row is exact by the 9-lemma.
Let $\bar y$ be the image of $y = \sigma(x)$ in $F_3/gF_3$, and let $\bar z$ be its image in $N$.
From the definition
of $N$, we see that  $\bar z$ is equal to the image
of $-\phi'(x)$, and is thus in $\gm N$.

Since $\bar y$ is a minimal generator of $F_3/gF_3$, we may factor out a complementary direct summand
from $F_3/gF_3$ and factor out its image in $N$ to get a short exact sequence of $S/(g)$-modules
$$
(*)\qquad 0\to (S/(g))\bar y\to N' \to I/(g) \to 0.
$$
 Let $\bar z'$ be the image of $\bar z$ in $N'$. Since $\bar z \in \gm N$, we have $\bar z' \in \gm N'$. 

We apply Theorem~\ref{gpit1} to $\bar z' \in N'$. Since $I$ has codimension 2, the module $I/(g)$ is generated by 1 element locally
at all the minimal primes of $(g)$, so, regarded as an $S/(g)$-module, $\rank(N') = 2$.
Thus there exists a prime ideal $P/(g)$ of $S/(g)$ that has codimension $\leq 2$ and contains $N'^*(\bar z') \subset S/(g)$. Since $P/(g)$ contains $N'^*(\bar z')$ the sequence $(*)$ does not split when localized at $P/(g)$.
In particular, $I/(g) \subset P/(g)$ so $I\subset P$.

We claim that $P$ contains the ideal $I(\phi_{3})$ generated by the maximal minors of $\phi_{3}$.
If this were not the case, the localization $(S/I)_P$ would be Cohen-Macaulay of codimension 2, and the module
$(I/(g))_P$, being the first syzygy of $(S/I)_P$ as an $(S/(g))_P$-module, would be a maximal Cohen-Macaulay $(S/(g))_P$-module.
Since $(S/(g))_P$ is Gorenstein, this would imply that the 
sequence (*) splits, a contradiction. 

Thus $I(\phi_3)\subset P$, and since $\codim P = 3$, Theorem~\ref{WMACE}, (3) implies that $P$ is an associated prime of $I$.
\end{proof}
\medbreak

\noindent\textbf{Conclusion of the proof of Theorem~\ref{main}: } We have now shown that if $I$ has codimension 2 and the projective dimension of $S/I$ is 3, then the failure of either condition in Criterion~\ref{criterion} leads to $I$ having
an embedded prime of codimension 3. This shows that if $I$ is unmixed then the ring $R = S/I$ is Golod, completing the proof of Theorem~\ref{main}.\qed

\bibliographystyle{amsalpha}
\bibliography{Golod}
\end{document}

%% file: preamble.tex
\def\antiddot{\mathinner{\mkern1mu\raise1pt\vbox{\kern7pt\hbox{.}}\mkern2mu
        \raise4pt\hbox{.}\mkern2mu\raise7pt\hbox{.}\mkern1mu}}

\newcommand{\BB}{{\mathbb B}}

\newcommand{\FF}{{\mathbb F}}

\newcommand{\PP}{{\mathbb P}}

\newcommand{\ZZ}{{\mathbb Z}}

\newcommand{\coker}{{\rm{coker}\,}}

\newcommand{\punkt}{\hspace{-.3ex}\raise.15ex\hbox to1ex{\Huge.}}

\def \fix#1 {{\hfill\break \bf (( #1 ))\hfill\break}}

\DeclareMathOperator{\syz}{syz}

\DeclareMathOperator{\depth}{depth}

\DeclareMathOperator{\pd}{pd}

\DeclareMathOperator{\Tor}{Tor}

\DeclareMathOperator{\codim}{codim}
\DeclareMathOperator{\rank}{rank}

\DeclareMathOperator{\height}{height}

\newcommand{\gm}{\mathfrak m}
\newcommand{\gn}{\mathfrak n}

\newtheorem{theorem}{Theorem}[section]

\newtheorem{corollary}[theorem]{Corollary}
\newtheorem{conjecture}[theorem]{Conjecture}
\theoremstyle{definition}
\newtheorem{definition}[theorem]{Definition}

\newtheorem{example}[theorem]{Example}